\documentstyle{amsart}

\numberwithin{equation}{section}

\newtheorem{theorem}{Theorem}[section]
\newtheorem{theorema}{Theorem}

\newtheorem{corollary}[theorem]{Corollary}
\newtheorem{proposition}[theorem]{Proposition}

\theoremstyle{remark}
\newtheorem{remark}[theorem]{Remark}   
          
\newtheorem{remarknon}{Remark}

\newtheorem{example}[theorem]{Example}
\newtheorem{ack}{Acknowledgment}

\newcommand{\ad}{\operatorname{ad}}

\newcommand{\Hom}{\operatorname{Hom}}

\newcommand{\g}{{\frak g}}

\begin{document}

\newcommand{\ra}{{\rho^\ast}}
\newcommand{\rb}{{\tilde{\rho}^\ast}}

\title{K\"ahler metrics on $G^{\Bbb C}$
}
\author{Roger Bielawski}
\thanks{Research supported by an EPSRC advanced fellowship}
%\dedicatory{preliminary version}

%\subjclass{53C25, 53C55}
\address{Department of Mathematics, University of Glasgow, Glasgow G12 8QW, UK }

\email{R.Bielawski@@maths.gla.ac.uk}

\begin{abstract} We study $G$-invariant K\"ahler metrics on  $G^{\Bbb C}$ from
the Hamiltonian point of view. As an application we show that there exist
$G\times G$-invariant Ricci-flat K\"ahler metrics on $G^{\Bbb C}$ for any
compact semisimple Lie group $G$.\end{abstract}

\maketitle

\section{Introduction}
We study $G$-invariant K\"ahler metrics on $G^{\Bbb C}$. We give Hamiltonian
ans\"atze for such metrics, in the spirit of \cite{LeB, PP} but for
non-abelian
 Lie groups. We are particularly interested in such metrics being Ricci-flat. We
give sufficient conditions for Ricci-flatness in terms of our Hamiltonian
 ansatz. For example, when $G$ is a Heisenberg group, we construct explicitly
$G$-invariant Ricci-flat K\"ahler metrics on a neighbourhood of $G$ in
$G^{\Bbb C}$. For compact $G$ we can do better:
%%%
\begin{theorema} Let $G$ be a compact simple Lie group and let $\gamma$ be a real closed $G\times G$-invariant $(1,1)$-form on $G^{\Bbb C}$. Then there exists a $G\times G$-invariant  K\"ahler
metric on $G^{\Bbb C}$ such that its Ricci form is $\gamma$. In particular,
there exists a $G\times G$-invariant Ricci-flat K\"ahler metric on $G^{\Bbb
C}$.\label{R-f-K}\end{theorema}
%%%
This result has been previously known for $G=SU(2)$ \cite{Can,Ste,Ste2,CGLP}.
In this case it is known that the Ricci-flat metric is complete and we expect
this to be true in the general case, although we have been unable to prove
it. Similarly the question of uniqueness is left open.
\par
The proof of Theorem \ref{R-f-K} relies on existence and regularity of {\bf
entire} solutions to the (real) Monge-Amp\`ere equation
%%%
\begin{equation} f(\nabla u)\det D_{ij}u=g\label{MA0}\end{equation}
where $f,g$ are certain smooth functions defined on ${\Bbb R}^n$. The ${\Bbb
R}^n$ in question is a Cartan subalgebra of $\text{Lie}(G)$ and all the
functions are invariant under the Weyl group $W$. There are very few results
about  entire solutions to \eqref{MA0}. In fact, it is usually assumed that
either $g\in L^1({\Bbb R}^n)$ or that $f=1$ and $g,1/g$ are bounded
\cite{CW}. None of these conditions applies in our situation and so we prove:
%%%%
\begin{theorema} Let $f,g:{\Bbb R}^n\rightarrow {\Bbb R}$ be two nonnegative
locally bounded functions invariant under a finite reflection group $W\subset
O(n)$ which acts irreducibly. If $\int_{{\Bbb R}^n}f(x)dx=+\infty$, then
there exists an entire convex $W$-invariant weak solution $u:{\Bbb
R}^n\rightarrow {\Bbb R}$ of \eqref{MA0}. If, in addition, $f$ and $g$ are
strictly positive and of class $C^{p,\alpha}$, then $u$ is of class
$C^{p+2,\alpha}$.\end{theorema}
%%%%
%%%%
Theorem \ref{R-f-K} can be viewed in a more general context. Given a
real-analytic Riemanian manifold $(M,g)$ one asks whether the metric $g$ can
be extended to a Ricci-flat K\"ahler metric $\bar{g}$ on a complex
thickening $M^{\Bbb C}$ of $M$, perhaps satisfying some additional
conditions. In this direction Bryant \cite{Bry} showed that if $\dim M=3$,
then such an extension exists with the original manifold $M$ being a special
Lagrangian submanifold. As Bryant indicates, such a result should also hold
for higher-dimensional manifolds with trivial tangent bundles. A similar
result is true when $(M,g)$ is K\"ahler \cite{Feix,Kal}. A stronger
condition on $\bar{g}$ is that there exists an anti-holomorphic and
isometric involution on $M^{\Bbb C}$ fixing $M$ (which is true in
\cite{Bry,Feix,Kal}). Our Theorem \ref{R-f-K}  (for $\gamma=0$) can be
viewed as an example of such an extension when $M=G$ with the bi-invariant
metric. The main point, however, is that our extension is {\it global}. In
this direction, one should ask a) when does there exist a geodesically
complete extension; b) what is the relation with the global existence of the
adapted global structure of Lempert and Sz\"oke \cite{LS} and of Guillemin
and Stenzel \cite{GS}. Recall that these authors have shown that a
real-analytic Riemannian manifold $(M,g)$ admits a canonical complex complex
structure, called {\it adapted}, on $TM$ or part thereof, characterised by
the condition that the geodesic foliation is holomorphic. Thus one could ask
whether there exists a Ricci-flat K\"ahler metric on the maximal domain of
definition of the adapted complex structure. The only previous cases where
we know of some sort global Ricci-flat K\"ahler extension are compact
homogeneous K\"ahler manifolds, where the Ricci-flat K\"ahler extensions are
actually hyperk\"ahler, and for compact symmetric spaces of rank $1$
\cite{Ste}.

\par
Although Theorem \ref{R-f-K} for $G=SU(2)$ has been proved by Candelas and de
la Ossa \cite{Can} and by Stenzel \cite{Ste}, our approach is somewhat
different and it expresses the K\"ahler potential directly as an invariant
function on $SL(2,{\Bbb C})$.
 %%%%
\begin{theorema} There exists a unique (up to homothety) complete $SU(2)\times SU(2)$
 invariant Ricci-flat K\"ahler metric on $SL(2,{\Bbb C})$. It is given by the K\"ahler
potential $$K(u)=\int_0^{R}\sqrt[3]{\sinh(2t)-2t}dt,$$ where $R=|h|$ and $h$
is an element of ${\frak su}(2)$ such that $u=g\exp{ih}$ with $g\in SU(2)$.
 \label{su(2)}  \end{theorema}

This metric is of cohomogeneity one. Ricci-flat K\"ahler metrics of
cohomogeneity one were studied by  Dancer and Wang \cite{DW} under the
assumption (generically satisfied) that the isotropy representation of the
principal orbit is multiplicity-free. The above metric is an example of
cohomogeneity one metric for which this assumption does not hold.

\section{Notation and conventions}
All $G$-actions are on the left. If $G$ acts on a smooth manifold $M$ and
$\rho$ is an element of the Lie algebra $\g$, then $\rho^\ast$ denotes the
fundamental vector field generated by $\rho$. As the action is on the left,
the map $\rho\mapsto\rho^\ast$ is an antihomomorphism (i.e.
$[\rho^\ast,\rb]=-[\rho,\tilde{\rho}]^\ast$).
\par
If $E$ is a representation of $G$, then by an $E$-valued $p$-form we mean a
$G$-equivariant map from $\Lambda^p(TM)$ to $E$, linear on fibers. If $E$ is
equipped with a $G$-invariant bilinear form $\langle\;,\;\rangle$ and
$\phi,\psi$ are $E$-valued $1$-forms on $M$, then we define a real-valued
$2$-form $\phi\wedge\psi$ by
\begin{equation} \phi\wedge\psi(X,Y)=\langle\phi(X),\psi(Y)\rangle -
\langle\phi(Y),\psi(X)\rangle.\end{equation} In general, we adopt the
convention that $v\wedge w=v\otimes w - w\otimes v$. Because of this and the
fact that we consider left $G$-actions, the structure equation for a
connection $1$-form $\theta$ on a principal $G$-bundle $P$ takes the form
\begin{equation} d\theta=[\theta,\theta]+\Omega\label{MC}\end{equation}
where $\Omega$ is the curvature of $\theta$.

\section{Complex structures on $G\times\g$}

Let $P$ be a (trivial) principal $G$-bundle over $\g$ and suppose that we are
given a $G$-invariant complex structure on $P$ (in other words we have a
$G$-invariant map $G^{\Bbb C}\rightarrow \g$). Any tangent vector on $P$ can
be written uniquely as $\ra+I\rb$ for some $\rho,\tilde{\rho}\in\g$. Since
$I$ is $G$-invariant, we have a connection $1$-form $\theta$ on $P$ defined
by
\begin{equation} \theta(\ra+I\rb)=\rho.\end{equation}
We also have a horizontal (i.e.  $\ad G$-valued and vanishing on vertical
vector fields) $1$-form $L$ defined by:
\begin{equation} L(\ra+I\rb)=\tilde{\rho}.\end{equation}
Conversely, given a connection $1$-form $\theta$ and a non-degenerate
horizontal
 $1$-form $L$ we can define an almost complex structure $I$ by:
\begin{equation} I\ra=\text{unique horizontal $Y$ such that
$L(Y)=\rho$}.\end{equation} We have
\begin{proposition} Let $P$ be a principal $G$--bundle over $\g$. A connection
$1$-form $\theta$ and a non-degenerate horizontal $1$-form $L$ on $P$ define
an integrable $G$-invariant complex structure on $P$ if and only if
\begin{equation} \begin{cases} \Omega=-[L,L]\\
DL=0.\end{cases}\label{SW}\end{equation}
\end{proposition}
\begin{remark} The minus sign is the consequence of $G$ acting on the left. For
$G$ semisimple, the assumption of horizontality of $L$ is unnecessary as it
follows from the first equation (and the non-degeneracy).\end{remark}
\begin{pf}
Since $I$ is $G$-invariant by definition, it satisfies $[\ra,IX]=I[\ra,X]$,
for any vector field $X$ and any $\ra\in\g$. It follows then from the formula
for the Nijenhuis tensor that $I$ is integrable if and only if
$[I\ra,I\rb]=[\rho,\tilde{\rho}]^\ast$ for any $\rho,\tilde{\rho}\in \g$.
From the definition of $I$ we have:
$$DL(I\ra,I\rb)=dL(I\ra,I\rb)=-L\bigl([I\ra,I\rb]\bigr),$$ where we have used
the fact, that for any $1$-form $\phi$ and any vector fields $X,Y$:
$$d\phi(X,Y)=X(\phi(Y))-Y(\phi(X))-\phi([X,Y]).$$
%%%
Therefore $[I\ra,I\rb]$ is vertical for all $\rho,\tilde{\rho}$ if and only
if $DL=0$. On the other hand, from the properties of curvature,
$$\Omega(I\ra,I\rb)=-\theta\bigl([I\ra,I\rb]\bigr)$$ and so the vertical part
of $[I\ra,I\rb]$ is equal to $[\rho,\tilde{\rho}]^\ast$ precisely when
$\Omega=-[L,L]$.
\end{pf}

As an example of a solution to \eqref{SW} consider a compact $G$ with the
$G$-equivariant diffeomorphism  between  $G\times \g$ and $G^{\Bbb C}$ given
by: $$ (g,h)\mapsto g\exp(ih),$$ i.e. by the polar decomposition. We shall
describe $\theta$ and $L$ given by this diffeomorphism. As everything is
equivariant it is enough to describe $\theta$ and $L$ at points of
$\{1\}\times \g$.
\par
For the exponential map of any Lie group we have the following formula:
$$(d\exp)_u\exp(-u)=\frac{\exp(\ad u)-1}{\ad u}.$$ Applying this to
$\exp{ih}$ and separating into the real and imaginary parts we obtain:
\begin{equation} \theta(\rho,v)=\rho+\frac{\cos(\ad h)-1}{\ad
h}(v),\label{theta}\end{equation}
\begin{equation} L(\rho,v)=\frac{\sin(\ad h)}{\ad h}(v)\label{L}\end{equation}
where $(\rho,v)$ is tangent to $G\times\g$ at the point $(1,h)$.

\section{K\"ahler metrics on $G^{\Bbb C}$}

We now consider the following problem. Let $\mu:G\times\g\rightarrow \g^\ast
$ be a regular $ G$-equivariant map. For a given
 $G$-invariant complex structure on $P=G\times\g$ defined by $\theta$ and $L$ we
wish to describe K\"ahler metrics on $P$ for which $\mu$ is the moment map.
\par
The covariant derivative $D\mu=d\mu+[\mu,\theta]$ vanishes on vertical vector
fields. This is also true for $L$ and since both $D\mu$ and $L$ are
non-degenerate, there exists an invertible $\Hom(\ad^\ast G,\ad G)$-valued
function $\Phi$ on $P$ such that
\begin{equation}L=\Phi(D\mu).\label{Phi}\end{equation}
A map $\Phi:\g^\ast\rightarrow \g$ can be viewed as giving a bilinear form on
$\g^\ast$: $\langle x, \Phi(y)\rangle$. We can therefore speak of $\Phi$
being symmetric, positive-definite etc. Let us introduce the following
notation. If $\mu$ is a map from a manifold into $\g^\ast$, then $\ad \mu$
denotes the map into $\Hom(\g, \g^\ast)$ defined by $\langle \ad
\mu(m)(x),y\rangle=\langle \mu(m),[x,y]\rangle$. Before stating the next
result, let us explain the notation used there. For any $m\in G\times\g$,
$\ad\mu\circ\Phi$ is a map from $\g^\ast$ to itself and we can talk about its
square. We have
%%%%
%%%%%
\begin{theorem} A regular equivariant map $\mu:P\rightarrow \g^\ast$ on
$P=G\times\g$ is a moment map for a $G$-invariant K\"ahler metric on $P$ if
and only if the function $\Phi$ defined by \eqref{Phi} is symmetric and if
both $\Phi$ and $\Phi +\Phi\circ(\ad\mu\circ\Phi)^2$ are positive-definite.
If this is the case then the K\"ahler metric is given by
\begin{equation}
g=\langle\Phi(D\mu),D\mu\rangle+\langle\Phi^{-1}(\theta),\theta\rangle +
\langle\mu,[\Phi(D\mu),\theta]-[\theta,\Phi(D\mu)]\rangle
\label{metric}\end{equation} and its K\"ahler form $\omega$ by
\begin{equation}
\omega=-d\langle\mu,\theta\rangle.\label{omega}\end{equation}
\end{theorem}
%%%%%
%%%%
\begin{pf} Suppose that $\mu$ is a moment map for a K\"ahler metric $g$. Then,
for any $\rho,\tilde{\rho}\in \g$,
\begin{equation} g(\ra,\rb)=\langle d\mu(I\ra),\tilde{\rho}\rangle= \langle
D\mu(I\ra),\tilde{\rho}\rangle=\langle\Phi^{-1}(\rho),\tilde{\rho}\rangle
\label{1}\end{equation} which shows that $\Phi$ is symmetric. Moreover
$$g(\ra,I\rb)=   \omega(\rb,\ra)=\langle d\mu(\ra),\tilde{\rho}\rangle=
\langle[\rho,\mu],\tilde{\rho}\rangle=
-\langle\mu,[\rho,\tilde{\rho}]\rangle$$ which shows that the metric has the
form \eqref{metric} (as $I\theta=-L$).
\par
To find out the conditions for this $g$ to be positive-definite we rewrite
$g$ as a metric on a Riemannian submersion. Let $\sigma$ be the connection
form of the Riemannian submersion defined by $g$. In other words
$\sigma(X)=\rho$ where $X-\ra$ is orthogonal to all vertical vector fields.
Since $\langle d\mu(IX), \tilde{\rho}\rangle=g(X,\rb)=g(\sigma(X)^\ast,\rb)$,
we have from \eqref{1}
\begin{equation} Id\mu=\Phi^{-1}(\sigma).\end{equation}
It follows that
\begin{equation} \sigma=I\Phi(d\mu)=I\Phi(D\mu-[\mu,\theta])=\theta
+\Phi\bigl([\mu,L]\bigr)\label{sigma}.\end{equation} Using this formula, we
can write metric \eqref{metric} as
%%%
\begin{equation} g=\langle\Phi(D\mu),D\mu\rangle -
\langle\Phi\bigl([\mu,\Phi(D\mu)]\bigr),[\mu,\Phi(D\mu)] \rangle +
\langle\Phi^{-1}(\sigma),\sigma\rangle. \label{submersion}\end{equation}
Computing the metric separately on vertical and horizontal (with respect to
$\sigma$) vectors, we see that both $\Phi$ and  $\Phi
+\Phi\circ(\ad\mu\circ\Phi)^2$ must be positive-definite.
\par
It remains to show that $\omega=g(I\;,\;)$ has the form \eqref{omega}. We
compute (using symmetry of $\Phi$ and \eqref{SW})
\begin{multline*} g(I\;,\;)= -\langle \Phi^{-1}(L),\theta\rangle + \langle
\Phi^{-1}(\theta),L\rangle +\langle\mu,[L,L]+[\theta,\theta]\rangle=\\
\theta\wedge D\mu+\langle\mu,-d\theta+2[\theta,\theta]\rangle = \theta\wedge
d\mu -\langle \mu,d\theta\rangle=-d\langle\mu,\theta\rangle.
\end{multline*}
Here $\phi\wedge\psi$ for a $\g$-valued $1$-form $\phi$ and a
$\g^\ast$-valued $1$-form $\psi$ denotes the $2$-form
$\phi\wedge\psi(X,Y)=\langle\psi(Y), \phi(X),\rangle -
\langle\phi(Y),\psi(X)\rangle$. This proves the theorem.
\end{pf}

If we are interested only in pseudo-K\"ahler metrics, then the relevant
condition is much simpler.
\begin{proposition} A regular equivariant map $\mu:P\rightarrow \g^\ast$ on
$P=G\times\g$ is a moment map for a $G$-invariant pseudo-K\"ahler metric on
$P$ if and only if the $1$-form $\langle\mu,L\rangle$ is
closed.\label{closed}
\end{proposition}
\begin{pf} From the proof of the above theorem, $\mu$ defines a pseudo-K\"ahler
metric if and only if $\Phi$ is symmetric ($\Phi$ is non-degenerate, since
$\mu$ is regular). Computing $d\langle\mu,L\rangle$ and using \eqref{SW}
shows that $d\langle\mu,L\rangle=D\mu\wedge L$. As $L=\Phi(D\mu)$, this last
expression vanishes precisely when $\Phi$ is symmetric.\end{pf}

\begin{example} Let $G$ be compact. There is a canonical $G\times G$ invariant
K\"ahler metric on $G^{\Bbb C}$ given by the K\"ahler potential
$K(u)=\frac{1}{2}|h|^2$, where $u=g\exp{ih}$, $g\in G$, $h\in \g$. The moment
map is given at points of $\{1\}\times\g$ by  $\mu(h)=h$. The K\"ahler form,
from  \eqref{omega} and  \eqref{theta}, is $$ \omega=-d\langle
h,\rho\rangle$$ where $\rho$ is the canonical flat connection on $G$. In
other words $\omega$ is just the canonical K\"ahler form of $T^\ast G=
G\times \g$. We also have $$D\mu=\cos(\ad h),$$ $$\Phi^{-1}=\cos(\ad
h)\frac{\ad h}{\sin(\ad h)}.$$ Using \eqref{submersion} one can easily
compute the metric on $G\backslash G^{\Bbb C}/G$, i.e. on a Weyl chamber,
which turns out to be the standard Euclidean metric.
\end{example}

\section{ Ricci-flat metrics and proof of Theorem 1}

As the Ricci curvature of a K\"ahler metric $g$ is given by
$-i\partial\bar{\partial} \ln \det g$, we wish to write the K\"ahler form in
terms of a holomorphic frame, namely
\begin{multline*} \omega= -d\mu\wedge\theta-\langle\mu,d\theta\rangle=
-\left(d\mu\wedge\theta+\langle\mu,[\theta,\theta]-[L,L]\rangle\right)=\\
-\frac{i}{2}\left(\Phi^{-1}(\theta+iL)\wedge
(\theta-iL)+i[\mu,\theta+iL]\wedge (\theta-iL)\right)=\\
-\frac{i}{2}\left((\Phi^{-1}+i\ad \mu)(\theta+iL)\right) \wedge
(\theta-iL).\end{multline*}
%%%%
Thus $g$ is Ricci-flat if the determinant of the operator $\Phi^{-1}+i\ad
\mu$ is constant for some basis of $\g$ and the dual basis.

\begin{example} We shall find invariant Ricci-flat K\"ahler metrics on the
complexification of a Heisenberg group. Let $(V,\omega)$ be a symplectic
vector space and let ${\frak h}(V)=V\times {\Bbb R}$ denote the corresponding
Heisenberg algebra. Thus $[(v,r),(w,s)]=(0,\omega(v,w))$. Let $H(V)$ denote
the corresponding (simply-connected) Lie group. To define a K\"ahler metric
on $H(V)^{\Bbb C}$ we need $\mu,\theta$ and $L$. Choose a symplectic basis
$p_1,q_1\dots,p_n,q_n$ of $V$, so that $$\omega=\sum p_i^\ast\wedge
q_i^\ast.$$ Then $p_1,q_1\dots,p_n,q_n,1$ is a basis of ${\frak h}(V)$ and we
define $\mu$ to be the map sending each vector of this basis to the
corresponding vector of the dual basis. We shall look for a $\Phi$ which, in
the above basis, is of the form $$\Phi_{(v,t)}=\text{diag}\;\bigl(1,\dots,
1,f(t)\bigr)$$ for a positive function $f$. It follows easily that
\begin{equation}\det (\Phi^{-1}+i\ad\mu)_{(v,t)}=\frac{(1-t^2)^n}{f(t)}.
\label{deter}\end{equation} The connection $1$-form $\theta$ will be given at
a point $(p_i,q_i,t)$ by $\bigl(0,0,\sum p_i dq^i\bigr)$. Then $D\mu=d\mu$
and $L=\Phi(d\mu)$ is given at a point $(p_i,q_i,t)$ by $(dp_i,dq_i,
f(t)dt)$. Since $[L,\theta]=0$ and $dL=0$, the equations \eqref{SW} are
satisfied. Thus we obtain an $H(V)$-invariant K\"ahler metric on $H(V)^{\Bbb
C}$ for any positive function $f(t)$. The equation \eqref{deter} implies that
the metric is Ricci-flat if $f(t)=c(1-t^2)^n$ for some constant $c$. These
metrics are defined only on an open subset of $H(V)^{\Bbb C}$ and are
incomplete there.
\end{example}

We shall now consider in detail the case of a compact semisimple $G$ and a
metric $g$ which is also invariant with respect to the right $G$-action.
First of all we have:
%%%
\begin{proposition} Let $G$ be compact semisimple and let $\gamma$ be a real closed $G\times G$-invariant $(1,1)$-form on $G^{\Bbb C}$. Then there exists a unique (up to an additive constant) real $G\times G$-invariant function $u$ such that $\gamma=-i\partial\bar{\partial}u$. \end{proposition}
%%%
\begin{pf} The existence of a function $u$ satisfying $\gamma=-i\partial\bar{\partial}u$
follows from the fact that $G^{\Bbb C}$ is Stein (and its second Betti
number is zero). As $G$ is compact we can average over $G\times G$ and
obtain a $G\times G$-invariant $u$. To show uniquness we have to prove that
a $G\times G$ invariant pluriharmonic ($\partial\bar{\partial}f=0$) function
$f$ on $G^{\Bbb C}$ is constant. We use the complex structure given by
$(\theta,L)$ in \eqref{theta} and \eqref{L}. Under this diffeomorphism
$G^{\Bbb C}\simeq G\times\g$, the right $G$-action on $G^{\Bbb C}$ becomes
the adjoint action on $\g$. Let $p_1,\dots,p_n$ be a basis of invariant
polynomials on $\g$. A $G\times G$ invariant function on $G\times \g$ can be
written as $f=f(p_1,\dots,p_n)$. We shall compute $dIdf$.
\par
Let us write $df=\langle F,dh\rangle$ for a map $F:\g\rightarrow \g$. Since
$dp_i([\rho,h])=0$ for any $\rho$ and at any point $h\in\g$, $[F(h),h]=0$ at
any regular $h$. Therefore, by \eqref{theta} and \eqref{L}, $Idf=\langle F,
dgg^{-1}\rangle$ and $$ dIdf= dF\wedge  dgg^{-1}+\langle F,[ dgg^{-1},
dgg^{-1}]\rangle.$$ Since $F$ is invariant for the left $G$-action,
$dF(\rho^\ast)$ vanishes for any $\rho\in \g$ and so $ \langle F,[ \rho,
\tilde{\rho}]\rangle=0$ for all $\rho, \tilde{\rho}\in\g$. Since $\g$ is
semisimple,  $F=0$ and $f$ is constant.
 \end{pf}

Therefore finding a K\"ahler metric with prescribed Ricci form is equivalent
to finding a K\"ahler metric of the form \eqref{metric} with prescribed
(positive) determinant of the hermitian operator $\Phi^{-1}+i\ad\mu$.

Now, notice that $$(\Phi^{-1}+i\ad \mu)\circ
L=D\mu+i[\mu,L]=d\mu+[\mu,\theta+iL].$$ Therefore, using $(\theta,L)$ given
by \eqref{theta} and \eqref{L} we require finding a map $\mu:\g\rightarrow
\g$ such that, at every point $h\in \g$, the determinant of the operator
\begin{equation} d\mu+\left[\mu, \frac{\exp(i\ad h)-1}{\ad h}\right]
\label{det}\end{equation} is equal to $e^u$-multiple of the determinant of
the operator \eqref{L}.
\par
 Since we look
for $G\times G$ invariant metrics on $G^{\Bbb C}$, $\mu$ must have a special
form. As the right action of $G$ becomes the adjoint action on $\g$ under the
polar decomposition and the two actions commute, the moment map $\mu$ for the
left action must be  $\ad G$-invariant. This means that $\mu$ maps any Cartan
subalgebra ${\frak h}$ to itself. Since, for any $x\in {\frak h}$ and any
$\rho\in \g$, $d\mu([\rho,x])=[\rho,\mu(x)]$, we easily compute the
determinant of \eqref{det} as $\bigl(\det d\mu_{|_{\frak h}}\bigr)\prod
\frac{\alpha (\mu(x))}{\alpha(x)}$, where the product is taken over all roots
$\alpha$. On the other hand, the operator $L$ given by $\eqref{L}$ has
determinant $\prod\frac{\sinh \alpha(x)}{\alpha(x)}$. Thus, to prove Theorem
1 we have to find a $W$-equivariant ($W$ being the Weyl group) map
$\mu:{\frak h}\rightarrow {\frak h}$ satisfying the equation
$$\left(\prod \alpha (\mu(x))\right)\det d\mu =e^{\tilde{u}(x)}\prod\sinh
\alpha(x)$$ for an arbitrary $W$-invariant function $\tilde{u}$.
Moreover, from Proposition \ref{closed}, $\langle\mu, L\rangle$ must be a
closed $1$-form on $\g$. As $\mu(x)$ commutes with $x$, this form is just
$\langle\mu(x), dx\rangle$. Therefore $\mu$ is the gradient of a
$G$-invariant function $K$ defined on $\g$. $K$ is of course the K\"ahler
potential of the metric. By restricting $K$ to ${\frak h}$ we obtain:
%%%
\begin{proposition}
For a compact semisimple Lie group $G$, the $G\times G$-invariant K\"ahler
metrics on $G^{\Bbb C}$ with Ricci form $\gamma=-i\partial\bar{\partial}u$
are (up to homothety) in one-to-one correspondence with smooth convex
$W$-invariant solutions $\tilde{K}$ to the Monge-Amp\`ere equation
\begin{equation}
\left(\prod \alpha (\nabla \tilde{K})\right)\det D_{ij}\tilde{K}=
e^{\tilde{u}(x)}\prod\sinh \alpha(x) \label{Winv}
\end{equation}
defined on an entire Cartan subalgebra ${\frak h}$, where the product is
taken over all roots $\alpha$ and $\tilde{u}=u_|{_{\frak h}}$.
\par
Equivalently, such metrics are in one-to-one correspondence with
 smooth convex $G$-invariant solutions $K$ to the Monge-Amp\`ere
equation
\begin{equation}
\det D_{ij}{K} = e^{u(x)}\frac{\prod\sinh \alpha(x)}{\prod
\alpha(x)}\label{Ginv}
\end{equation}
defined on all of $\g$ (where the right-hand side is viewed as an invariant
function on $\g$). \label{yeah}\end{proposition}
%%%%
%%%%%
\begin{pf}
We have to show that $\Phi$ and $\Phi\circ (\ad \mu\circ \Phi)^2$ are
positive definite. As $D\mu$ is given by the real part of \eqref{det}, it
follows that $L$ and $D\mu$ commute. As they are both self-adjoint, so is
$\Phi=L\circ(D\mu)^{-1}$. $L$ is positive definite and $D\mu$ is positive
definite on vectors tangent to ${\frak h}$, where it is given by the matrix
of second derivatives of $\tilde{K}$. On the other hand $D\mu$ is diagonal on
the root spaces and its eigenvalues are $\frac{ \alpha (\nabla
\tilde{K})}{\alpha(x)}\cosh\alpha(x)$. Since $\tilde{K}$ is convex, these
eigenvalues are positive. Therefore $\Phi$ is positive definite. To check
that  $\Phi +\Phi\circ(\ad\mu\circ\Phi)^2$ is positive-definite note that
$\Phi$ and $\ad \mu$ commute. Therefore we have $$ \Phi
+\Phi\circ(\ad\mu\circ\Phi)^2= \Phi^3(\Phi^{-1}+i\ad\mu)
(\Phi^{-1}-i\ad\mu)$$ which is positive definite as  $\Phi^{-1}+i\ad\mu$ is
hermitian (it is nondegenerate as $\det (\Phi^{-1}+i\ad\mu)$ is
positive).\end{pf}
%%%%%
%%%%
\begin{remark} For $u=0$ (i.e. the Ricci-flat case), the Monge-Amp\`ere equation \eqref{Ginv} has the following
interpretation. We seek a map $\mu:\g\rightarrow \g$ with a convex potential
(i.e. $\mu=\nabla K$ for a convex function $K$) such that the pullback under
$\mu$ of the canonical volume form of $\g$ (given by the Killing metric) is
equal to the volume form $\hat{\omega}$ which is the pullback of the
canonical volume form on the symmetric space $G\backslash G^{\Bbb C}$ via the
exponential map $\exp:\g\rightarrow G\backslash G^{\Bbb C}$.\end{remark}

\begin{corollary}
Let $G$ be a compact simple Lie group and let $\gamma$ be a real closed
$G\times G$-invariant $(1,1)$-form on $G^{\Bbb C}$. Then there exists a
$G\times G$-invariant K\"ahler metric on $G^{\Bbb C}$ whose Ricci form is
$\gamma$.
\end{corollary}
\begin{pf}
 Theorem \ref{Mat} in the
appendix shows the existence of a weak solution $\tilde{K}$ to \eqref{Winv}.
By Proposition \ref{proper}, $\tilde{K}$ is proper. Since a weak solution is
equivalent to a solution a.e. (every convex function has first and second
derivatives a.e.), $\tilde{K}$ gives rise to a proper (convex) weak solution
${K}$ of \eqref{Ginv} on $\g$. To show that ${K}$ is smooth apply Theorem 2
in \cite{Cafinterior} to $u(x)=K -c$, $c\in {\Bbb R}$ and the convex set
$\Omega_c=\{x; K(x)\leq c\}$.\end{pf}

We have been unable to show that the Ricci-flat K\"ahler metrics obtained
from Theorem 1 are complete. From \eqref{submersion}, it follows easily
\begin{proposition} Let $G$ be a compact semisimple Lie group of rank $n$.
 A $G\times G$-invariant
K\"ahler metric on $G^{\Bbb C}$ given by a $W$-invariant convex solution
$\tilde{K}$ to \eqref{Winv} is complete if and only if
$$\sum\frac{\partial^2\tilde{K}}{\partial x_i\partial x_j} dx_i\otimes dx_j$$
is a complete metric on ${\Bbb R}^n$. \hfill $\Box$\end{proposition}

\section{ Ricci-flat K\"ahler metrics on $SU(2)$}

We will compute explicitly $SU(2)\times SU(2)$-invariant metrics on
$SL(2,{\Bbb C})$.  By the arguments of the previous section such a metric is
given by a moment map of the form  $\mu(h)=f(|h|)h$ for some real function
$f$ (at points of  $\{1\}\times {\frak su}(2) $). We compute the determinant
of the operator \eqref{det}. It can be rewritten as
\begin{equation} f(|h|)1 +\dot{f}(|h|)d(|h|)h +f(|h|)(\cos(\ad h)-1+i\sin(\ad
h)).\label{operator}\end{equation} Now, on ${\frak su}(2)$, $(\ad h)^2$ has
the eigenvalue $0$ corresponding to the eigenvector $h$ and the eigenvalue
$-|h|^2$ with the eigenspace corresponding to vectors orthogonal to $h$ (we
fix here an $\ad SU(2)$-invariant metric on ${\frak su}(2)$, so that the norm
of the Pauli matrices is $1$). For the operator \eqref{operator} $h$ is also
an eigenvector and its eigenvalue is: $$\lambda_0=f(|h|)+\dot{f}(|h|)|h|.$$
We can also easily compute the other two eigenvalues as being
$$\lambda_{\pm}=f(|h|)\bigl(\cosh |h|\pm \sinh|h|\bigr)= f(|h|)e^{\pm|h|}.$$
On the other hand, the operator $L$ given by $\eqref{L}$ has eigenvalues $1$
and
 (twice) $\sinh(|h|)/|h|$. Therefore the metric will be Ricci flat if and only
if the function $f=f(t)$ satisfies
%%%%
\begin{equation}  f^2(f+t\dot{f})=c\left(\frac{\sinh t}{t}\right)^2
\label{equation} \end{equation}
%%%%
where $c$ is a  positive constant. The constant $c$ amounts to changing the
metric by a homothety and so we can assume that $c=1$. Let us denote the
function on the right by $G(t)$ and let us also write $u=f^3/3$. Then the
last equation becomes $$t\frac{du}{dt}+3u=G(t),$$ whose general solution is
$$ u(t)=\left(a+\int_0^t s^2G(s)ds\right)t^{-3}. $$
%%%%%%%
%%%%
As $u$ needs to be smooth at the origin, $a=0$. Therefore
\begin{equation} u(t)=t^{-3}\int_0^t(\sinh\tau)^2d\tau=
\frac{\sinh(2t)-2t}{4t^3} \label{solution}\end{equation}

To check that we obtained a K\"ahler metric we have to, according to Theorem
3.1, check that $\Phi$ is self-adjoint and that both $\Phi$ and $\Phi
+\Phi\circ(\ad\mu\circ\Phi)^2$ are positive-definite. This is done in the
proof of Proposition \ref{yeah}.
\par
Thus we have obtained an $SU(2)\times SU(2)$ invariant Ricci-flat K\"ahler
metric on $SL(2,{\Bbb C})$. It follows from the arguments that any other such
a metric is homothetic to this one.
\par
Let us show that this metric is complete. As it is  $SU(2)\times SU(2)$
invariant, it is of cohomogeneity one and it is enough to show that the
metric on the quotient ${\Bbb R}_{\geq 0}$ is complete. This means computing
the first two terms in \eqref{submersion} in the radial direction $h/|h|$.
The second term becomes zero. The first one is $\langle L,D\mu\rangle$, and
since $L(h/|h|)=h/|h|$ and
$D\mu(h/|h|)=\bigl(f(|h|)+\dot{f}(|h|)|h|\bigr)h/|h|$ (from \eqref{operator}
and the subsequent paragraph), the quotient metric is just
$$(t\dot{f}+f)dt^2.$$
%%%%
From \eqref{equation} $t\dot{f}+f= G(t)(3u)^{-2/3}$ and it follows from
\eqref{solution} that $t\dot{f}+f$ has an exponential growth. Therefore the
metric is complete.
\par
Finally, let us compute its K\"ahler potential, i.e. a real-valued function
$K$ such that $\omega=-i\partial\bar{\partial} K$. By Theorem 3.1,
$\omega=-d\langle\mu,\theta\rangle$. Therefore a K\"ahler potential will be a
function $K$ such that $2\langle\mu,\theta\rangle=IdK$. i.e. $dK=2\langle\mu,
L\rangle$. In our case the right-hand side becomes $ 2\langle
f(|h|)h,dh\rangle$, which means that $K=K(|h|)$ with $K(t)$ satisfying
$\frac{dK}{dt}=tf$. This shows that the K\"ahler potential has the form
stated in Theorem \ref{su(2)} (up to a constant multiple).

\appendix\section{Entire $W$-invariant solutions of Monge-Amp\`ere equations}

We wish to show existence and regularity of entire solutions to a class of
Monge-Amp\`ere equations
\begin{equation} f(\nabla \phi)\det D_{ij}\phi=g(x)\label{MA} \end{equation}
where $f$ and $g$ are nonnegative functions on ${\Bbb R}^n$. We recall first
the concept of a weak solution of \eqref{MA}. Let $\phi$ be a convex
function. Then $\nabla \phi$ is a well-defined multi-valued mapping:
$(\nabla\phi)(x)$ is the set of slopes of all supporting hyperplanes to the
graph of $\phi$ at $(x,\phi(x))$. If $B$ is a subset of ${\Bbb R}^n$, let
$\nabla\phi(B)$ be its image in the multi-valued sense. Then $\phi$ is a weak
solution of \eqref{MA} if $$\int_B g(x)dx=\int_{\nabla\phi(B)}f(y)dy$$ for
every Borel set $B$. Let us denote the right-hand side by $\omega(B,\phi,f)$.
It can be shown that it is a Borel measure on ${\Bbb R}^n$. A basic result is
that if $u_k\rightarrow u$ compactly and $f_k \rightarrow f$ uniformly, then
$\omega(\cdot,u_k,f_k)$ converges to $\omega(\cdot,u,f)$ weakly, i.e. as
functionals on the space of compactly supported continuous functions.
\par
After these preliminaries we are going to prove:
%%%%%
\begin{theorem} Let $W\subset O(n)$ be a finite reflection group acting
irreducibly on ${\Bbb R}^n$ and let $f,g$ be two nonnegative $W$-invariant
functions on ${\Bbb R}^n$. Furthermore, assume that $f$ and $g$ are locally
bounded and that $$\int_{{\Bbb R}^n}f=+\infty.$$ Then there exists a (weak)
$W$-invariant solution $\phi:{\Bbb R}^n\rightarrow {\Bbb R}$ of the
Monge-Amp\`ere equation \eqref{MA}. Moreover $\phi$ is convex and Lipschitz
continuous.\label{Mat}
\end{theorem}
%%%%
%%%%
\begin{pf}
Put $f_k=f+1/k$ and $g_k=g+1/k$, $k\in {\Bbb N}_+$. Let $B_r$ denote the ball
of radius $r$ centred at the origin. Let $R_k$ be a number defined by
$$\int_{B_{R_k}}f_k=\int_{B_k}g_k.$$ According to \cite{Bre,Cafconvex} there
exists a unique (up to a constant) strictly convex solution $\phi_k$ of
\begin{equation} f_k(\nabla \phi_k)\det D_{ij}\phi_k=g_k(x)\label{MAk}
\end{equation}
which is of class $C^{1,\beta}$, for some $\beta>0$, and such that $\nabla
\phi_k$ maps $B_k$ onto $B_{R_k}$. Moreover, as $W\subset O(n)$ and all the
data is $W$-invariant, $\phi_k$ is $W$-invariant (this follows from
uniqueness, since $\phi_k\circ w$ is also a solution of \eqref{MAk}). To
prove the existence of a weak solution to \eqref{MA}, defined on all of
${\Bbb R}^n$, it is enough to show that the functions $\phi_k$ are uniformly
(in $k$) bounded on any ball $B_R$, $k\geq R$. Indeed,
 a bounded sequence of convex functions on a bounded open convex domain has
a convergent subsequence. This follows from the elementary fact, which we
will use repeatedly, that the slopes of supporting hyperplanes of a convex
function, bounded by $R$ on a domain $G$, are bounded by $2R/\delta$ on any
subdomain $G^\prime$ such that $\text{dist}(G^\prime,\partial G)\geq \delta$.
\par
Let us show that the functions $\phi_k$, $k\geq R$, are bounded on $B_R$
uniformly in $k$. Let $\Delta$ be a minimal set of reflections generating $W$
and consider a subgroup $W^\prime$ of $W$ generated by $n-1$ elements of
$\Delta$. Let $L$ be the one-dimensional subspace of ${\Bbb R}^n$ fixed by
$W^\prime$. Since $\phi_k$ is $W$-invariant, $\nabla \phi_k$ maps $L$ to
itself. We claim that $\nabla \phi_k$ are uniformly bounded on $L\cap B_R$.
Suppose that they are not. Then there exists a subsequence $k_j$ such that
$\nabla \phi_{k_j}$ maps $L\cap B_R$ onto $L\cap B_{s_j}$ ($\phi_k$, being
$W$-invariant, satisfy $\nabla\phi_k(0)=0$), where $s_j\rightarrow +\infty$.
Now consider the unique solution $\psi_k$, $\psi_k(0)=0$, to the equation $$
g_k(\nabla \psi_k)\det D_{ij}\psi_k=f_k(y)$$ mapping $B_{R_{k_j}}$ onto
$B_{k_j}$. According to \cite{Bre}, $\nabla \psi_k$ is the inverse of $\nabla
\phi_k$. Thus $\nabla \psi_{k_j}$ maps  $L\cap B_{s_j}$ onto $L\cap B_R$. It
follows that the functions $\psi_{k_j}$ are uniformly bounded (by $rR$) on
any $L\cap B_r$ (for sufficiently large $j$). As the $\psi_{k_j}$ are
$W$-invariant and convex, they are bounded by $rR$ on the convex hull of the
set $W(L\cap B_r)$. Since $W$ acts irreducibly, $WL$ generates ${\Bbb R}^n$
and so the sets  $\text{conv}\left(W(L\cap B_r)\right)$ cover ${\Bbb R}^n$.
In fact, there is a number $\delta\in (0,1)$, such that $B_{\delta r}\subset
\text{conv}\left(W(L\cap B_r)\right)$ for any $r$. It follows that there is a
subsequence of $\psi_{k_j}$ convergent to a convex solution $\psi:{\Bbb
R}^n\rightarrow {\Bbb R}$ of $$g(\nabla\psi)\det D_{ij}\psi=f.$$ The function
$\psi$ is bounded by $Rr/\delta$ on any $B_r$. Since $\psi$ is convex,
$\nabla\psi$ is bounded by $4R/\delta$ on $B_{r/2}$. Thus
$U=\nabla\psi\bigl({\Bbb R}^n\bigr)$ is contained in $B_{4R/\delta}$. However
$$+\infty=\int_{{\Bbb R}^n}f= \int_U g, $$ which leads to a contradiction.
\par
We have shown so far that $\nabla\phi_k$ are uniformly bounded on $L\cap
B_R$. Therefore the $\phi_k$ are bounded on $L\cap B_R$ (we assume
$\phi_k(0)=0$). The argument applied before to $\psi_k$ can now be used for
the $\phi_k$: since $WL$ generates ${\Bbb R}^n$ and $\phi_k$ are convex and
$W$-invariant, $\phi_k$ are bounded on every compact subset of ${\Bbb R}^n$.
Thus there exists a Lipschitz continuous limit of some subsequence.
\end{pf}

\begin{proposition}  Suppose, in addition, that $\int_{{\Bbb R}^n}g>0$. Then
any entire $W$-invariant convex solution $\phi$ of \eqref{MA}  is a proper
function. \label{proper}
\end{proposition}
\begin{pf}
 First of all, the assumption on $g$ implies that  $\phi$ cannot
be a bounded function. Indeed, otherwise $\phi$,  being convex and defined
 on all of ${\Bbb R}^n$, is constant.  \newline
 Let $C_W$ denote the polytopal complex of $W$, i.e. all hyperplanes fixed by
 reflections in $W$ and all their intersections.
We proceed by induction on $\dim L$, $L\in C_W\cup \{{\Bbb R}^n\}$, to show
that $\phi$ is proper on all elements of $C_W\cup \{{\Bbb R}^n\}$.
 Let $L$ be a $1$-dimensional element of $C_W$.
 Then $\phi$ is proper on $L$,
 since otherwise $\phi$ would be bounded (as $\text{conv}\left(WL)\right)$
 is all of ${\Bbb R}^n$). Now assume that $\phi$ is proper on all
 $p$-dimensional elements of $C_W$. Let $L\in C_W$, $\dim L=p+1$. We claim
 that $\phi$ is actually proper on the vector space $V_L$ spanned by $L$.
Indeed,
 otherwise there is a number $K$ and a set of points $x_k\rightarrow\infty$
 such that $\phi(x_k)\leq K$. For sufficiently large $k$, the points $x_k$ do not
 lie on $p$-dimensional
 elements of $C_W$ contained in $L$. There is a reflection subgroup
 $W^\prime\subset W$ acting on $V_L$. By taking the intervals joining points
 of $W^\prime x_k$, for each $k$, we obtain a contradiction, as these
 intervals intersect the $p$-dimensional
 elements of $C_W$ contained in $L$ and $\phi$ is convex.
\end{pf}

\begin{corollary}
In the situation of Theorem \ref{Mat} suppose, in addition,  that $f$ and $g$
are strictly positive and of class $C^{p,\alpha}$. Then any entire
$W$-invariant convex solution $\phi$ of \eqref{MA} is
$C^{p+2,\alpha}$.\label{smooth}
\end{corollary}
\begin{remarknon} It is well-known that the conclusion of this corollary does not
hold for solutions which are not $W$-invariant.\end{remarknon}
\begin{pf}
We can apply Corollary 2 in \cite{Caflocalization} to the convex (and
bounded by the last proposition) set $\Omega_c=\{x;\phi(x)\leq c\}$. Doing
this for every $c$ shows that $\phi$ is strictly convex. Now the main result
of \cite{Cafsome} implies that $\phi$ is of class $C^{1,\beta}$. Therefore
$\nabla \phi$ is H\"older continuous and $\phi$ is a solution of $\det
D_{ij}\phi=\tilde g$, where $\tilde{g}=g/f(\nabla \phi)$ is of class
$C^{0,\beta}$. We can  apply Theorem 2 in \cite{Cafinterior} to
$u(x)=\phi(x) -c$, $c\in {\Bbb R}$ and the convex and bounded set
$\Omega_c=\{x; u(x)\leq 0\}$ to conclude that $\phi$ is of class
$C^{2,\beta}$. Higher regularity is standard.
\end{pf}

\begin{ack} It is a pleasure to thank David Calderbank for comments; in
particular for his suggestion that the Hamiltonian ansatz of section 4 can be
stated in a form suitable for non-compact groups.\end{ack}


\begin{thebibliography}{99}

\bibitem{Bre}
{Y. Brenier} `Polar factorization and monotone rearrangement of vector-valued
functions', {\it Comm. Pure Appl. Math.} 44 (1991),  375--417.

\bibitem{Bry}
{R.L. Bryant} `Calibrated embeddings in the special Lagrangian and
coassociative cases', {\it Ann. Global Anal. Geom.} 18 (2000), 405--435.

\bibitem{Caflocalization}
{L.A. Caffarelli} `A localization property of viscosity solutions to the
Monge-Amp\`ere equation and their strict convexity', {\it Ann. of Math. (2)}
131 (1990), 129--134.


\bibitem{Cafinterior}
{L.A. Caffarelli} `Interior $W\sp {2,p}$ estimates for solutions of the
Monge-Amp\`ere equation', {\it Ann. of Math. (2)} 131 (1990), 135--150.


\bibitem{Cafsome}
{L.A. Caffarelli} `Some regularity properties of solutions of Monge-Amp\`ere
equation', {\it Comm. Pure Appl. Math.} 44 (1991), 965--969.


\bibitem{Cafconvex}
{L.A. Caffarelli} `The regularity of mappings with a convex potential',{\it
J. Amer. Math. Soc.} 5 (1992), 99--104.

\bibitem{Can}
{P. Candelas \and X.C. de la Ossa} `Comments on conifolds', {\it Nuclear
Phys. B} 342 (1990), 246--268.

\bibitem{CW}
{K.-S. Chou \and X.-J. Wang} `Entire solutions of the Monge-Amp\`ere
equation', {\it Comm. Pure Appl. Math.} 49 (1996),  529--539.

\bibitem{CGLP}
 {M. Cvetic, G.W. Gibbons, H. Lu \and C.N. Pope}
 `Ricci-flat metrics, harmonic Forms and brane resolutions', hep-th/0012011.

\bibitem{DW}
{A. Dancer \and M. Wang}, `K\"ahler-Einstein metrics of cohomogeneity one',
{\it Math. Ann.} 312 (1998),  503--526.

\bibitem{Feix}
{B. Feix}, `Hyperk\"ahler metrics on cotangent bundles', {\it J. Reine Angew.
Math.} 532 (2001), 33--46.


\bibitem{GS}
{V. Guillemin \and M.B. Stenzel} `Grauert tubes and the homogeneous
Monge-Amp\`ere equation', {\it J. Differential Geom.} 34 (1991), 561--570.


\bibitem{Kal}
{D. Kaledin}, `A canonical hyperk\"ahler metric on the total space of a
cotangent bundle', in `Quaternionic structures in mathematics and physics
(Rome, 1999)', 195--230, World Sci. Publishing, River Edge, NJ, 2001.

\bibitem{LeB}
{C. LeBrun}, `Explicit self-dual metrics on ${\Bbb C}P^2\#\cdots\#{\Bbb
C}P^2$',
 {\it J. Differential Geom.}
34 (1991), 223--253.


\bibitem{Ste2}
{T.-C. Lee}, `Complete Ricci flat K\"ahler metric on $M\sp n\sb {\rm I}$,
$M\sp {2n}\sb {\rm II}$, $M\sp {4n}\sb {\rm III}$', {\it Pacific J. Math.}
185 (1998),  315--326.

\bibitem{LS}
{L. Lempert \and R. Sz\"oke} `Global solutions of the homogeneous complex
Monge-Amp\`ere equation and complex structures on the tangent bundle of
Riemannian manifolds', {\it Math. Ann.} 290 (1991), 689--712.

\bibitem{PP}
{H. Pedersen \and Y.S. Poon}, `Hamiltonian constructions of K\"ahler-Einstein
metrics and K\"ahler metrics of constant scalar curvature', {\it Comm. Math.
Phys.} 136 (1991), 309--326.


\bibitem{Ste}
{M.B. Stenzel} `Ricci-flat metrics on the complexification of a compact rank
one symmetric space', {\it Manuscripta Math.} 80 (1993),  151--163.





\end{thebibliography}
\end{document}